# A FLEMING–VIOT PROCESS AND BAYESIAN NONPARAMETRICS

BY STEPHEN G. WALKER, SPYRIDON J. HATJISPYROS
AND THEODOROS NICOLERIS

*University of Kent, University of the Aegean and University of the Aegean*

This paper provides a construction of a Fleming–Viot measure valued diffusion process, for which the transition function is known, by extending recent ideas of the Gibbs sampler based Markov processes. In particular, we concentrate on the Chapman–Kolmogorov consistency conditions which allows a simple derivation of such a Fleming–Viot process, once a key and apparently new combinatorial result for Pólya-urn sequences has been established.

**1. Introduction.** The Fleming–Viot process, introduced by Fleming and Viot [6], is a measure valued diffusion process. The stationary distribution of the process is $\Pi$, where $\Pi$ is the distribution of a random measure $\mu$, on some space $S$, and $\mu$ can be obtained via

$$\mu(\cdot) = \sum_{i=1}^{\infty} \rho_i \, \delta_{V_i}(\cdot), \tag{1.1}$$

where $\rho_1 > \rho_2 > \cdots$ have the Poisson–Dirichlet distribution [8] and $V_1, V_2, \ldots$ are independent and identically distributed from $\nu_0$, and independent of the $\rho_i$. Such a random measure is also known as a Dirichlet process [5] and has been of great importance to Bayesian nonparametric methods. To denote the dependence on $(\theta, \nu_0)$, we will use the notation $\Pi(\theta \nu_0)$.

Ethier and Griffiths [4] provide the transition function for a particular Fleming–Viot process. Let $d_n(t) = P(D_t = n)$, where $D_t$ is a death process, $D_0 = \infty$ a.s., and with rate $\lambda_n = \frac{1}{2}n(n-1+\theta)$ for some $\theta > 0$. Tavaré [11], for example, computed that, for $n = 1, 2, \ldots,$

$$d_n(t) = \sum_{m=n}^{\infty} (-1)^{m-n} C(m,n)(\theta+n)_{(m-1)} m!^{-1} \gamma_{m,t,\theta}, \tag{1.2}$$

where

$$\gamma_{m,t,\theta} = (2m - 1 + \theta)e^{-\lambda_m t}$$

and

$$d_0(t) = 1 - \sum_{m=1}^{\infty} (-1)^{m-1} \theta_{(m-1)} m!^{-1} \gamma_{m,t,\theta}.$$

Also,

$$C(m,n) = \frac{m!}{(m-n)!n!}$$

and $a_{(m)} = a(a+1)\cdots(a+m-1)$ for $m = 1, 2, \ldots$ with $a_{(0)} = 1$. We will also use $a_{[m]} = a(a-1)\cdots(a-m+1)$ for $m = 1, 2, \ldots$ with $a_{[0]} = 1$. We will show among other things that this death process is fundamentally connected with the general Pólya-urn scheme [3].

The transition function is given by

$$(1.3) \quad P(t, \mu, d\nu) = \sum_{n=0}^{\infty} d_n(t) \int \Pi\left(d\nu \Big| \theta\nu_0 + \sum_{i=1}^{n} \delta_{X_i}\right) \mu(dX_1) \cdots \mu(dX_n).$$

It is the intention of this paper to establish a comprehensive construction of the process and the transition function using ideas formulated in Bayesian nonparametrics relating to the Dirichlet process. The key result, which appears new and involves an elegant combinatorial identity, is for sequences of Pólya-urns. We will also use recent ideas for constructing Markov processes using latent variables, outlined in [10].

In Section 2 we provide background to the construction of the Fleming–Viot process via discrete time processes associated with the Dirichlet process. The necessary Chapman–Kolmogorov condition for existence in continuous time is examined in Section 3. Section 4 contains technical results and Section 5 concludes the paper with some points of discussion.

**2. Stationary Markov processes using the Dirichlet process.** In [10] the use of the Dirichlet process for deriving the DAR(1) model was described. Consider the joint distribution on $S \times \mathscr{P}(S)$, where $\mathscr{P}(S)$ is the space of probability measures on $S$, given by

$$P(d\mu, dX) = \mu(dX)\Pi(d\mu|\theta\nu_0).$$

In words, $\mu$ is chosen from $\Pi$ and, given $\mu$, $X$ is chosen from $\mu$. By making use of both conditional distributions, the conditional distribution for $\mu$ being

$$P(d\mu|X) = \Pi(d\mu|\theta\nu_0 + \delta_X),$$

a discrete time Markov process can be constructed on $S$ with transition function

$$P(X_t, dX_{t+1}) = \int \mu(dX_{t+1})\Pi(d\mu|\theta\nu_0 + \delta_{X_t}).$$

This is of the form

$$P(X_t, dX_{t+1}) = \int P(dX_{t+1}|\mu)P(d\mu|X_t).$$

A result in [2] gives

(2.1) $$\Pi(\theta\nu_0) = \int \Pi(\cdot|\theta\nu_0 + \delta_X)\nu_0(dX)$$

and it is well known that

(2.2) $$\int \mu\Pi(d\mu|\theta\nu_0) = \nu_0.$$

Consequently, it is easy to show that $\nu_0$ is the stationary distribution of the process. The process, using properties of the Gibbs sampler, is easily shown to be reversible.

In fact, it follows that

$$X_{t+1} \begin{cases} \sim \nu_0, & \text{with probability } \theta/(1+\theta), \\ = X_t, & \text{with probability } 1/(1+\theta). \end{cases}$$

This is the DAR(1) model.

Alternatively, we could consider the measure valued process on $\mathscr{P}(S)$. This would have transition function given by

$$P(\mu, d\nu) = \int \Pi(d\nu|\theta\nu_0 + \delta_X)\mu(dX).$$

Using (2.1) and (2.2), it is straightforward to show that $\Pi$ is the stationary distribution of the process and that it is also reversible.

Instead of having a single observation from $\mu$, we could consider the joint distribution on $S^n \times \mathscr{P}(S)$ given by

$$P(dX_1, \ldots, dX_n, d\mu) = \Pi(d\mu|\theta\nu_0)\prod_{i=1}^n \mu(dX_i).$$

It is well known that the "posterior" or conditional distribution of $\mu$ given $X_1, \ldots, X_n$ has the form

$$\Pi\left(\cdot\bigg|\theta\nu_0 + \sum_{i=1}^n \delta_{X_i}\right);$$

see [5]. Hence, in this case, the transition function for the measure valued process is

$$(2.3) \qquad P(\mu, d\nu) = \int \Pi\bigg(d\nu \bigg| \theta\nu_0 + \sum_{i=1}^{n} \delta_{X_i}\bigg) \mu(dX_1) \cdots \mu(dX_n).$$

This is beginning to resemble the transition function for the Fleming–Viot process given in (1.3), though obviously for discrete time. To obtain the Fleming–Viot process, we need to put the processes we have constructed into continuous time. To achieve this, it is necessary to make the number of samples $n$ be random. So, using the notation of Ethier and Griffiths in [4], we denote by $d_n(t)$ the probability that the number of samples being used is $n$ for the transition with time $t$.

We are now ready to examine the existence of such a process by looking at the Chapman–Kolmogorov equations.

**3. Chapman–Kolmogorov conditions.** We consider the transition from $\mu_0$ to $\mu_t$ and then to $\mu_{t+s}$. Thus, we have

$$P(d\mu_{t+s}|Y_1,\ldots,Y_m,m) = \Pi\bigg(d\mu_{t+s}\bigg|\theta\nu_0 + \sum_{j=1}^{m}\delta_{Y_j}\bigg),$$

where $Y_1,\ldots,Y_m$ are independent and identically distributed from $\mu_t$ and $P(m) = d_m(s)$. Also,

$$P(d\mu_t|X_1,\ldots,X_n,n) = \Pi\bigg(d\mu_t\bigg|\theta\nu_0 + \sum_{i=1}^{n}\delta_{X_i}\bigg),$$

where $X_1,\ldots,X_n$ are independent and identically distributed from $\mu_0$ and $P(n) = d_n(t)$. Now, it is well known that we can integrate out $\mu_t$ and that

$$P(Y_1,\ldots,Y_m|X_1,\ldots,X_n,n) = \mathscr{Q}\bigg(\theta\nu_0 + \sum_{i=1}^{n}\delta_{X_i}\bigg),$$

where $\mathscr{Q}$ denotes a distribution associated with the general Pólya-urn scheme [2, 3]. In terms of sampling, we take the $Y_1, Y_2, \ldots$ by sampling $Y_1 \sim \nu_n$, where

$$\nu_n = \frac{\theta\nu_0 + \sum_{i=1}^{n}\delta_{X_i}}{\theta + n},$$

and then, subsequently,

$$P(Y_j|Y_1,\ldots,Y_{j-1}) = \frac{(\theta+n)\nu_n + \sum_{l=1}^{j-1}\delta_{Y_l}}{\theta + n + j - 1}.$$

Such a sampling scheme, which in the Bayesian nonparametric literature is known as the Pólya-urn scheme, also appears in mathematical genetics and is connected with the Poisson–Dirichlet process. See, for example, [11]. Since the $X_i$ are independent and identically distributed from $\mu_0$, to achieve the Chapman–Kolmogorov condition, we need to understand how many of the $Y_j$'s are independent and identically distributed from $\mu_0$. For there to be $r$ of them, we obviously require $m \geq r$ and $n \geq r$. To expand on this, a number of the $Y_j$'s will be identical to some of the $X_i$'s. We are looking for the number of distinct indices associated with these $X_i$'s. So, for example, if we collect up $\{X_2, X_4, X_2, X_1, X_6, X_4\}$ from sampling the $Y_j$'s, then the appropriate number is $r = 4$; that is, we have the distinct indices $\{1, 2, 4, 6\}$.

THEOREM 3.1. *Conditionally on $n$ and $m$, we have that the probability mass function for $r$ is given, for $r \in \{0, \ldots, \min\{n, m\}\}$, by*

$$P(r|m,n) = \frac{n_{[r]}(\theta + r)_{(m-r)}}{(\theta + n)_{(m)}} C(m,r),$$

*which can be written in extended form as*

$$P(r|m,n) = r! C(n,r) C(m,r) \frac{\theta_{(n)} \theta_{(m)}}{\theta_{(n+m)} \theta_{(r)}}.$$

The proof is provided in Section 4. Hence, the Chapman–Kolmogorov condition becomes

$$(3.1) \qquad d_r(t+s) = \sum_{n=r}^{\infty} \sum_{m=r}^{\infty} r! C(n,r) C(m,r) \frac{\theta_{(n)} \theta_{(m)}}{\theta_{(n+m)} \theta_{(r)}} d_m(s) d_n(t)$$

for all $s, t > 0$. There can be many solutions to this, we will look for those within the class of death processes; that is, $d_n(t) = P(D_t = n)$. We let the rate be $\lambda_n$ and so, in particular, we have $P(D_{s+t} = n | D_t = n) = P(T_n > s)$, where $T_n$ is an exponential r.v. with parameter $\lambda_n$.

We have the death process also satisfying Chapman–Kolmogorov and so

$$d_r(t+s) = \sum_{n=r}^{\infty} P(D_{t+s} = r | D_t = n) d_n(t).$$

Comparing this with the Chapman–Kolmogorov condition in (3.1) for the measure valued process, we see that we should have

$$\frac{\theta_{(r)} P(D_{t+s} = r | D_t = n)}{r!} = \sum_{m=r}^{\infty} C(m,r) \frac{\theta_{(m)}}{\theta_{(n+m)}} d_m(s) \theta_{(n)} C(n,r)$$

and so

$$(3.2) \qquad \sum_{m=r}^{\infty} C(m,r) \frac{\theta_{(m)}}{\theta_{(n+m)}} d_m(s) = \frac{\theta_{(r)}(n-r)!}{\theta_{(n)} n!} P(D_{t+s} = r | D_t = n).$$

This needs to be solved.

We will now show that the $d_n(t)$ given in (1.2) is a solution to (3.2). Ethier and Griffiths [4] did much the same thing, but from a more complicated starting point. Our demonstration, which follows, is now straightforward given the result of Theorem 3.1.

Now

(3.3) $$P(D_{t+h} = n | D_t = n) = 1 - \lambda_n h + o(h)$$

and

(3.4) $$P(D_{t+h} = n-1 | D_t = n) = \lambda_n h + o(h).$$

By considering

$$\sum_{m=n}^{\infty} \frac{m_{[n]}}{(\theta+m)_{(n)}} d_m(s)$$

and

$$\sum_{m=n-1}^{\infty} \frac{m_{[n-1]}}{(\theta+m)_{(n)}} d_m(s)$$

which are part of (3.2) with $r = n$ and $r = n-1$, respectively, with the help of formulae appearing in [4], page 1585, we can show that the $d_m(s)$ given in (1.2) satisfies the conditions for the death process. The details are provided in Section 4.

Hence, we see that the complicated nature of the death process probabilities is solely due to the form of $P(r|n,m)$, which is a property of the Pólya-urn scheme. For other processes, perhaps with a different choice of $\Pi$, which generates discrete random distribution functions (see Section 5.2), and so yield different $P(r|n,m)$, the fundamental equation to solve for obtaining a transition function satisfying Chapman–Kolmogorov is to find $d_m(s)$ such that

$$P(D_{t+s} = r | D_t = n) = \sum_{m=r}^{\infty} P(r|n,m) d_m(s).$$

This appears to be the key.

**4. Technical results.**

*Result* [A]. We first prove that

$$\sum_{m=n}^{\infty} \frac{m_{[n]}}{(\theta+m)_{(n)}} d_m(s) = e^{-\lambda_n s}.$$

Now the left-hand side can be written as

$$\sum_{m=n}^{\infty}\sum_{k=m}^{\infty}(-1)^{k-m}C(k,m)m_{[n]}k!^{-1}\frac{(\theta+m)_{(k-1)}}{(\theta+m)_{(n)}}\gamma_{k,s,\theta},$$

which is equal to

(4.1)
$$\frac{n_{[n]}\gamma_{n,s,\theta}}{n!(\theta+2n-1)}$$
$$+\sum_{k=n+1}^{\infty}\frac{\gamma_{k,s,\theta}}{(k-n)!}\sum_{m=n}^{k}(-1)^{k-m}C(k-n,m-n)(\theta+m+n)_{(k-n-1)}.$$

LEMMA 4.1. *It is that*

$$\sum_{l=0}^{k}(-1)^{k-l}C(k,l)(\phi+l)_{(k-1)}=0$$

*for any $\phi>0$ and $k\geq 1$.*

PROOF. We will do this by induction and prove the more general result that

$$\sum_{l=0}^{k}(-1)^{k-l}C(k,l)(\phi+l)_{(k-r)}=0$$

for all $r\in\{1,\ldots,k\}$. Assume the result is true for all $k<K$ and for $r\in\{R,\ldots,K\}$ when $k=K$. Now

$$\sum_{l=0}^{K}(-1)^{K-l}C(K,l)(\phi+l)_{(K-R+1)}$$

$$=\phi\sum_{l=0}^{K}(-1)^{K-l}C(K,l)(1+\phi+l)_{(K-R)}$$

$$+K\sum_{l=1}^{K}(-1)^{K-l}C(K-1,l-1)(1+\phi+l)_{(K-1-R+1)},$$

which by hypothesis is zero. To complete the proof, note that

$$\sum_{l=0}^{k}(-1)^{k-l}C(k,l)=0$$

for all $k=1,2,\ldots$ and that the result is true for $K=2$. □

The result follows from (4.1) by setting $l=m-n$, and then substituting $k-n$ for $k$. Hence, $d_m(\cdot)$ satisfies (3.3).

*Result* [B]. We next show, in the first instance, that

$$H(s) = \sum_{m=n-1}^{\infty} \frac{m_{[n-1]}}{(\theta+m)_{(n)}} d_m(s) = \frac{1}{2} \frac{1}{\lambda_n - \lambda_{n-1}} (e^{-\lambda_{n-1}s} - e^{-\lambda_n s}).$$

Using a result in [4], page 1585, we have

$$\frac{dH(s)}{ds} = \frac{1}{2} e^{-\lambda_{n-1}s} - \lambda_n H(s)$$

and so

$$H(s) = \frac{1}{2} \frac{1}{\lambda_n - \lambda_{n-1}} e^{-\lambda_{n-1}s} + C e^{-\lambda_n s}.$$

Now $H(0) = 0$ and so

$$C = -\frac{1}{2} \frac{1}{\lambda_n - \lambda_{n-1}},$$

leading to the desired result. Now performing some elementary algebra on (3.2) with $r = n-1$, we obtain

$$P(D_{t+s} = n-1 | D_t = n) = n(n-1+\theta) H(s),$$

which leads to the validity of (3.4).

We could have proven the validity of (3.3) using this technique as well. If

$$G(s) = \sum_{m=n}^{\infty} \frac{m_{[n]}}{(\theta+m)_{(n)}} d_m(s),$$

then

$$\frac{dG(s)}{ds} = -\lambda_n G(s).$$

Hence, $G(s) = C \exp(-\lambda_n s)$ and since $G(0) = 1$, we have the result.

*Result* [C]. Before proving Theorem 3.1, we need to establish the following result (this is apparently a new combinatorial result):

LEMMA 4.2. *Let $\theta$ be a positive real. Then for $m \geq r > 0$ and defining $0! = 1$,*

$$\sum_{k=0}^{m-r} k! C(k+r-1,k) C(m-r,k) \theta_{(m-r-k)} = (\theta+r)_{(m-r)}.$$

PROOF. Let $|s(n,k)|$ denote the *unsigned* or *absolute* Stirling numbers of the first kind. Expanding the $\theta$ terms on both sides of this relation, we obtain

$$\sum_{k=0}^{m-r} k! C(k+r-1,k) C(m-r,k) \sum_{l=0}^{m-r-k} |s(m-r-k,l)| \theta^l$$

$$= \sum_{k=0}^{m-r} |s(m-r,k)| \sum_{l=0}^{k} C(k,l) \theta^l r^{k-l}.$$

By changing the order of summation on both sides and collecting up the terms, we have

$$\sum_{l=0}^{m-r} \left\{ \sum_{k=0}^{m-r-l} k! C(k+r-1,k) C(m-r,k) |s(m-r-k,l)| \right\} \theta^l$$

$$= \sum_{l=0}^{m-r} \left\{ \sum_{k=l}^{m-r} C(k,l) |s(m-r,k)| r^{k-l} \right\} \theta^l.$$

These are two polynomials in $\theta$ of degree $m-r$ and for them to be equal $\forall \theta$, it suffices to establish the equality of the coefficients of the same powers of $\theta$; that is,

$$\sum_{k=0}^{m-r-l} k! C(k+r-1,k) C(m-r,k) |s(m-r-k,l)| = \sum_{k=l}^{m-r} C(k,l) |s(m-r,k)| r^{k-l},$$

for all $l = 0, 1, \ldots, m-r$.

To show this, is true, we make use of an identity that appears in [1], page 824, which states that, for positive integers $a \leq b \leq c$,

$$C(b,a) |s(c,b)| = \sum_{j=b-a}^{c-a} C(c,j) |s(c-j,a)| |s(j,b-a)|.$$

The right-hand side of the equation can be written as

$$\sum_{k=l}^{m-r} C(k,l) |s(m-r,k)| r^{k-l}$$

$$= \sum_{k=l}^{m-r} \sum_{j=k-l}^{m-r-l} C(m-r,j) |s(m-r-j,l)| |s(j,k-l)| r^{k-l}$$

$$= \sum_{j=0}^{m-r-l} C(m-r,j) |s(m-r-j,l)| \sum_{k=l}^{j+l} |s(j,k-l)| r^{k-l}$$

$$= \sum_{j=0}^{m-r-l} C(m-r,j)|s(m-r-j,l)| \sum_{k=0}^{j} |s(j,k)| r^k$$

$$= \sum_{j=0}^{m-r-l} C(m-r,j)|s(m-r-j,l)| r_{(j)}$$

$$= \sum_{j=0}^{m-r-l} j! C(j+r-1,j) C(m-r,j) |s(m-r-j,l)|,$$

where we have used the fact that $r_{(j)} = j! C(j+r-1,j)$. This completes the proof. □

PROOF OF THEOREM 3.1. The closed form of the joint probability of
$$[Y_1, \ldots, Y_m | X_1, \ldots, X_n]$$
is given by

$$dG(Y_1, \ldots, Y_m | X_1, \ldots, X_n) = \prod_{i=1}^{m} \left\{ \frac{\theta \nu_0 + \sum_{j=1}^{i-1} \delta_{Y_j}(dY_i) + \sum_{l=1}^{n} \delta_{X_l}(dY_i)}{\theta + n + i - 1} \right\}.$$

Assume without loss of generality that $X_1, \ldots, X_r$, the first $r$ observations from $X_1, \ldots, X_n$, are those that are repeated when we obtain a sample $Y_1, \ldots, Y_m$. Let $0 \le s_i \le m - r, i = 1, 2, \ldots, r$, and fix a number $k$ such that $s_1 + s_2 + \cdots + s_r = k$, where $0 \le k \le m - r$.

Here the $s_i$ represent the multiplicity of the $X_i$, $i = 1, 2, \ldots, r$, that appear in the sample when there are $k$ spaces available for those repetitions. So, conditionally on $X_1, \ldots, X_n$, we are searching for the probability of the simultaneous occurrence of the following events:

$$Y_1 = X_1, \qquad Y_2 = X_2, \ldots, Y_r = X_r,$$
$$Y_{r+j} = X_1, \qquad 1 \le j \le s_1,$$
$$Y_{r+s_1+j} = X_2, \qquad 1 \le j \le s_2,$$
$$\vdots$$
$$Y_{r+\sum_{i=1}^{r-1} s_i + j} = X_r, \qquad 1 \le j \le s_r,$$
$$Y_{r+k+1} \in \mathcal{X} - \{X_1, \ldots, X_n\},$$
$$Y_{r+k+j} \in \mathcal{X} - \{Y_{r+k+1}, \ldots, Y_{r+k+j-1}, X_1, \ldots, X_n\} \quad \text{or}$$
$$Y_{r+k+j} \in \{Y_{r+k+1}, \ldots, Y_{r+k+j-1}\}, \qquad 2 \le j \le m - r - k,$$

where $\mathcal{X}$ is the sample space. This probability is given by
$$\frac{(s_1+1)! \cdots (s_r+1)! \theta_{(m-r-k)}}{(\theta+n)_{(m)}}.$$

Since these events are exchangeable, by taking into consideration the number of repetitions of the $X_i$'s and the specific order of appearance of the new values, which depend on the previous observations, then, for a given $k$ and for fixed multiplicities $s_1, \ldots, s_r$, the probability of them occurring in any order is given by

$$\left\{\frac{m!}{(s_1+1)!\cdots(s_r+1)!(m-r-k)!}\right\}\frac{(s_1+1)!\cdots(s_r+1)!\theta_{(m-r-k)}}{(\theta+n)_{(m)}}$$

$$= \frac{m!\theta_{(m-r-k)}}{(m-r-k)!(\theta+n)_{(m)}}.$$

If we let $k$ and $s_1, \ldots, s_r$ vary, then this probability becomes

$$\sum_{k=0}^{m-r} \sum_{\{s_1+\cdots+s_r=k\}} \frac{m!\theta_{(m-r-k)}}{(m-r-k)!(\theta+n)_{(m)}}$$

$$= \sum_{k=0}^{m-r} C(k+r-1,k) \frac{m!\theta_{(m-r-k)}}{(m-r-k)!(\theta+n)_{(m)}}$$

$$= \frac{m!}{(m-r)!} \sum_{k=0}^{m-r} k! C(k+r-1,k) C(m-r,k) \frac{\theta_{(m-r-k)}}{(\theta+n)_{(m)}}$$

and since, as proven in Lemma 4.2,

$$\sum_{k=0}^{m-r} k! C(k+r-1,k) C(m-r,k) \theta_{(m-r-k)} = (\theta+r)_{(m-r)},$$

this probability becomes

$$\frac{m!(\theta+r)_{(m-r)}}{(m-r)!(\theta+n)_{(m)}}.$$

Finally, for any choice of $r$ $X$'s from $\{X_1, \ldots, X_n\}$, we have that

$$P(r|m,n) = C(n,r)\frac{m!(\theta+r)_{(m-r)}}{(m-r)!(\theta+n)_{(m)}},$$

which is given by

$$\frac{n_{[r]}(\theta+r)_{(m-r)}}{(\theta+n)_{(m)}} C(m,r),$$

as required. □

**5. Discussion.** We have shown how to construct a particular Fleming–Viot process, for which the transition function is known, from basic ideas involving the Dirichlet process and Markov processes, based on the Gibbs sampler. This approach requires a new combinatorial result involving Pólya-urn schemes. In particular, the combinatorial complexities which arise with the generator approach are avoided with the Chapman–Kolmogorov condition, once Theorem 3.1 has been established. Here we briefly discuss a number of points:

5.1. *The case $\theta = 0$.* Here we consider the case when $\theta = 0$. It is evident that since $\Pi$ no longer exists in this case, there can be no stationary distribution for the process. A stationary distribution, which is $\Pi(\theta \nu_0)$, can only exist when $\theta > 0$. When $\theta = 0$, the death process has probabilities

$$d_n(t) = \sum_{m=n}^{\infty} (-1)^{m-n} C(m,n) n_{(m-1)} m!^{-1} \gamma_{m,t},$$

for $n \geq 2$, where $\gamma_{m,t} = (2m-1) \exp\{-m(m-1)t/2\}$, with $d_0(t) = 0$ and

$$d_1(t) = 1 - \sum_{m=2}^{\infty} (-1)^m \gamma(m,t).$$

Now $d_n(t)$ is the probability that there are $n$ equivalence classes at time $t$ in the coalescent of [9]. When $\theta = 0$, then

$$P(Y_1, \ldots, Y_m | X_1, \ldots, X_n, n) = \mathscr{Q}\left(\sum_{i=1}^{n} \delta_{X_i}\right)$$

and

$$P(r|n,m) = \frac{n_{[r]} r_{(m-r)}}{n_{(m)}} C(m,r),$$

which can be written as

$$P(r|n,m) = rC(m,r)C(n,r)\frac{(n-1)!(m-1)!}{(n+m-1)!},$$
$$n = 1, 2, \ldots; m = 1, 2, \ldots.$$

Hence, for $n, m > 0$, we have $P(r = 0 | n, m) = 0$.

5.2. *The next step.* We believe the representation given is informative, making a strong connection between Bayesian nonparametrics and population genetics. It is also based on first principles for the construction of a Markov process, namely, the proposal for a transition function and the verification of the Chapman–Kolmogorov condition. What are the possible directions in which this connection can be taken? The clear idea is that we can consider alternative choices of $\Pi$ which generates discrete random distribution functions. One class of such a random distribution function can be generated via

$$\mu(\cdot) = \sum_{i=1}^{\infty} \rho_i \delta_{V_i}(\cdot),$$

where the $V_i$ are independent and identically distributed from some measure $\nu_0$ and the $\rho_i$ have a stick-breaking structure; that is,

$$\rho_1 = w_1 \quad \text{and} \quad \rho_i = w_i \prod_{j<i}(1 - w_j),$$

where the $w_j$ have independent beta distributions, say, $\text{beta}(\alpha_j, \beta_j)$. Then $\mu$ is almost surely a random probability measure when

$$\sum_{j=1}^{\infty} \log(1 + \alpha_j/\beta_j) = \infty;$$

see [7]. For example, the Dirichlet process arises when $\alpha_j = 1$ and $\beta_j = \theta$. The two parameter Poisson–Dirichlet process, which is worth exploring, arises when $\alpha_j = 1 - \sigma$ and $\beta_j = \theta + j\sigma$ for $0 < \sigma < 1$ and $\theta > -\sigma$. To find the transition function for this process and others, if they exist, we would need to replicate Theorem 3.1, that is, find the appropriate $P(r|m, n)$ from the predictive distributions, and then solve

$$P(D_{t+s} = r | D_t = n) = \sum_{m=r}^{\infty} P(r|n, m) d_m(s)$$

for an appropriate death process. Hence, we have a strategy for finding alternative transition functions which seems to be highly possible to achieve. Work on this is ongoing.

5.3. *An inequality.* Here we consider the usefulness of

$$\sum_{m=n}^{\infty} \frac{m_{[n]}}{(\theta + m)_{(n)}} d_m(t) = e^{-\lambda_n t}.$$

For example, by putting $n = 1$, we have

$$\sum_{m=1}^{\infty} \frac{m}{\theta + m} d_m(t) = e^{-\lambda_1 t}.$$

Hence, it is easy to obtain

$$e^{-\lambda_1 t} < 1 - d_0(t) < (1 + \theta) e^{-\lambda_1 t}$$

and it is also clear how to obtain improved inequalities from this identity.

**Acknowledgments.** The first author is grateful to Dario Spano for introducing him to the paper of Ethier and Griffiths [4]. The authors are grateful to two referees and the Editor for suggesting ways of improving the paper. The first author is supported by an EPSRC Advanced Research Fellowship.

S. G. WALKER  
INSTITUTE OF MATHEMATICS  
STATISTICS AND ACTUARIAL SCIENCE  
UNIVERSITY OF KENT  
CT2 7NZ CANTERBURY  
UNITED KINGDOM  
E-MAIL: S.G.Walker@kent.ac.uk

S. J. HATJISPYROS  
T. NICOLERIS  
DEPARTMENT OF STATISTICS AND  
  ACTUARIAL–FINANCIAL MATHEMATICS  
UNIVERSITY OF THE AEGEAN  
GR-832 00 KARLOVASSI, SAMOS  
GREECE  
E-MAIL: schatz@aegean.gr  
       nikthe@aegean.gr